\documentclass[]{amsart}

\usepackage{amsmath, amssymb,amsthm}
\usepackage{amsaddr}

\newtheorem{lemma}{Lemma}[section]
\newtheorem{theorem}[lemma]{Theorem}

\newtheorem{proposition}[lemma]{Proposition}

\theoremstyle{definition} 
\newtheorem{definition}[lemma]{Definition}
\newtheorem{remark}[lemma]{Remark}
\newtheorem{example}[lemma]{Example}

\DeclareMathOperator{\InDer}{InDer}
\DeclareMathOperator{\InStr}{InStr}
\DeclareMathOperator{\spin}{spin}
\DeclareMathOperator{\End}{End}
\DeclareMathOperator{\Rep}{\mathbf{Rep}}
\DeclareMathOperator{\ad}{ad}

\title[Semisimplifying Lie algebras of $J$-ternary algebras]{Semisimplifying Lie algebras of $J$-ternary algebras in characteristic $3$}
\author{Michiel Smet}
\email{Michiel.Smet@UGent.be}
\address{Department of Electronics and Information Systems\\ Faculty of Engineering and
	Architecture, Ghent University}



\begin{document}

\begin{abstract}
We describe a class of Lie superalgebras in characteristic $3$, containing the Elduque-Cunha superalgebras $\mathfrak{g}(3,3), \mathfrak{g}(6,6)$ and the Elduque superalgebra $\mathfrak{el}(5,3)$, using the tensor product of composition algebras. 
For the Lie superalgebra $\mathfrak{el}(5,3)$, this allows us to move beyond the contragredient construction and it also allows us to construct more general forms. 

We also describe how one obtains these Lie superalgebras using the semisimplification functor on the representation category $\mathbf{Rep}(\alpha_3)$ to Lie algebras of type $E_6, E_7$ and $E_8$, in line with how Arun Kannan applied this functor to the split algebras. We further apply this functor more broadly to the class of Lie algebras coming from $J$-ternary algebras over fields of characteristic $3$. 
\end{abstract}

\maketitle

\section*{Introduction}

Over fields of low characteristic, one can encounter simple Lie superalgebras for which there exist no analogs in other characteristics. In \cite{CUN07, Eld08}, the Freudenthal magic square was extended to incorporate composition superalgebras (which exist only over fields of characteristic $3$), yielding the first construction of the Lie superalgebras $\mathfrak{g}(3,3)$ and $\mathfrak{g}(6,6)$.  Of the simple Lie superalgebras, the ones with indecomposable Cartan matrix have been classified over fields with arbitrary characteristics \cite{Bou06,Bou09}, leading to the discovery of the Lie superalgebra $\mathfrak{el}(5,3)$.
We construct a family of simple Lie superalgebras over fields of characteristic $3$ containing $\mathfrak{g}(3,3), \mathfrak{g}(6,6), \mathfrak{el}(5,3)$ in the notation of \cite{Bou09}. 
We introduce a uniform construction using the tensor product of composition algebras. In particular, this yields a construction of $\mathfrak{el}(5,3)$ which is not contragredient and further allows for non-split forms of this Lie algebra to be constructed.

A more abstract construction than the one we develop (which will be explained in some detail as motivation), can be achieved using the semisimplification of certain categories \cite{KAN21}. 
The algebras already mentioned correspond to the images under semisimplification of $\mathfrak{e_6},\mathfrak{e_7}, \mathfrak{e_8}$. This relies on Allison's construction of isotropic Lie algebras of type $E_6$, $E_7$ and $E_8$ as $\mathbb{Z}$-graded Lie algebras with only weights $\pm 2, \pm 1, 0$ using the tensor product of composition algebras  \cite[8.4]{ALL79}. We apply semisimplification to the representation of $\alpha_3$ corresponding to an element $v$ in the Lie algebra which has weight $\pm 2$ such that $v$ is invertible in the structurable algebra.
In the case that the composition algebras are split, this agrees with parts of what Kannan \cite{KAN21} did with these Lie algebras, as the Lie superalgebras one obtains are contragredient. 
The broader idea to obtain superalgebras over a field of low characteristic from algebras using the semisimplification functor has been more widely applied to obtain supercomposition algebras and a super variant of Albert algebras, see e.g. \cite{Eld24, Eld25}.

We also apply the construction in a more general setting. Hein and Allison \cite{HEIN75, ALL76} introduced $J$-ternary algebras to study $\mathbb{Z}$-graded Lie algebras $L = J_{-2} \oplus M_{-1} \oplus D_0 \oplus M_1 \oplus J_2$ with $J_{\pm 2}$ copies of a Jordan algebra $J$ and $M_{\pm 1}$ copies of a $J$-module. These algebras are easily seen to be $\textbf{Rep}(\alpha_3)$-algebras using $k \mapsto [k 1_2, \cdot]$, where $1_2 \in J_2$ is the unit of the Jordan algebra. For simple $J$-ternary algebras, which are nearly classified over fields of characteristic different from $2$ (one broad class remains not fully classified when the characteristic is $3$), we carry out the construction of the Lie superalgebra. In this light, we also interpret the results of \cite{ELD06} as applying the semisimplification functor to the class of $J$-ternary algebras of degree $1$.\\

Very recently, Cunha and Elduque \cite{cunha2025j} independently had the same idea: describing the semisimplification of the Lie algebra associated to a $J$-ternary algebra, and using the tensor product of composition algebras to describe the algebras $\mathfrak{g}(3,3), \mathfrak{g}(6,6), \mathfrak{el}(5,3)$. In \cite[Theorem 5.1]{cunha2025j}, they describe what Lie superalgebras one obtains from the ``prototypical" $J$-ternary algebras, which we do in 1.4.1 and 1.4.2, using a more concrete way. By generalizing results from \cite{All88}, they can semisimplify using slightly more general elements. 
\\

Throughout, we assume that $\Phi$ is a field of characteristic different from $2$.
By the category of $\Phi$-algebras $\Phi\textbf{-alg}$ we mean the category of commutative, unital, associative $\Phi$-algebras.

\section{Semisimplifying $J$-ternary algebras}

In this section, we explain (i) how one semisimplifies the category $\Rep(\alpha_3)$, (ii) what $J$-ternary algebras are, (iii) the relation between $J$-ternary algebras of a certain form, (iv) what the semisimplification of the associated Lie algebra looks like, and (v) what the results of semisimplification looks like in the easy cases.

\subsection{Semisimplification of $\text{Rep}(\alpha_3)$}

Here, we give an application of the theory of semisimplification specifically applied to $\textbf{Rep}(\alpha_3)$. For a more general approach, see, e.g., \cite{Ostrik2015OnSF, etingof2021semisimplification} or \cite[section 3]{KAN21}.

Consider the functor $\alpha_3 : \Phi\textbf{-alg} \longrightarrow \textbf{Grp}$ defined by \[\alpha_3(K) = \{ x \in K | x^3 = 0\}\] with group operation $(x,y) \mapsto x + y$ over a field $\Phi$ of characteristic $3$. A representation of such a group functor is a vector space $V$ endowed with actions
\[ \alpha_3(K) \times (K \otimes V) \longrightarrow K \otimes V\]
which are natural in $K$, i.e., if $f : K \longrightarrow L$ is a homomorphism,
\[ (f \otimes 1)(g \cdot (k \otimes v)) = \alpha_3(f)(g) \cdot (f(k) \otimes v)\]
for all $k \in K$, $v \in V$ and $g \in \alpha_3(K)$. 

\begin{lemma}
	The $\alpha_3$ representations on a vector space $V$ are in a one-to-one correspondence with linear maps $f : V \longrightarrow V$ such that $f^3 = 0$.
	\begin{proof}
		Using the action of $\alpha_3(\Phi[t]/(t^3))$ on $V$ we obtain
		\[ t \cdot v = 1 \otimes v_0 + t \otimes v_1 + t^2 \otimes v_2,\]
		with $v \mapsto v_i$ linear in $v$.
		The naturality of the action implies that
		\[ g \cdot v = 1 \otimes v_0 + g \otimes v_1 + g^2 \otimes v_2\]
		for all $g \in \alpha_3(K)$ and all $K$.
	   We now define $f(v) = v_1$ for $v \in V$. Using that $(s + t) \cdot v = s \cdot (t \cdot v)$ for the action of $\alpha_3(\Phi[s,t]/(s^3,t^3))$ we see that $f(v_1) = 2 v_2, f(v_2) = 0$ and thus $f^3(v) = 0$.
	   Furthermore, we know that $v_0 = 0 \cdot v = v$. Hence, we conclude
	   \[ g \cdot v = v + g \otimes f(v) + g^2 \otimes f(f(v))/2.\]
	   
	   Conversely,
	   \[ k \cdot v = v + k f(v) + k^2 f(f(v))/2\]
	   defines an $\alpha_3$-representation for all linear $f$ such that $f^3 = 0$.
	\end{proof}
\end{lemma}

\begin{lemma}
	Suppose that $V$ is an $\alpha_3$-representation with respect to an $f$ such that $f^3 = 0$.
	A linear map $m : V \otimes V \longrightarrow V : a \otimes b \longrightarrow ab$ is $\alpha_3$-invariant if and only if 
	\[ f(ab) = f(a)b + af(b).\]
	Hence $(V,m)$ is an $\alpha_3$-algebra with respect to $f$ if and only if $f^3 = 0$ and $f$ is a derivation.
	\begin{proof}
		This is a straightforward computation.
	\end{proof}
\end{lemma}

 The indecomposable objects of $\Rep(\alpha_3)$ are given by $\Phi^a$ for $a = 1$, $2$, and $3,$ each with basis $\{e_i | i \le a\}$ on which the $f$ is given by $f(e_i) = e_{i+1}$, writing $e_l = 0$ for $l > a$.
 The only simple object, however, is $\Phi^1$. Hence the category $\textbf{Rep}(\alpha_3)$ is not semisimple. 
We can semisimplify this category \cite[Proposition 2.3]{Ostrik2015OnSF} to obtain a semisimple category $\overline{\textbf{Rep}(\alpha_3)}$ with as simple objects $\overline{\Phi}, \overline{\Phi^2}$. 
This category is obtained by dividing out all \textit{negligible} morphisms, i.e., morphisms $g:X \longrightarrow Y$ such that $\text{Tr}(gu) = 0$ for all $u : Y \longrightarrow X$.
We remark that $\overline{\Phi^3} \cong 0$ since all endomorphisms of $\overline{\Phi^3}$ have trace $0$ and only the object $0$ can have the $0$-morphism as identity.

\begin{remark}
	Given that we work with specific tensor categories, identifying what the key structural data are, is more important than a full definition. The important data are the Hom-spaces, which should be $\Phi$-vectorspaces and the tensor product $\otimes$. A braiding on a tensor category is a natural isomorphism $c_{X,Y} : X \otimes Y \longrightarrow Y \otimes X$ satisfying certain naturality conditions and axioms. For precise definitions of tensor categories and braidings, see \cite[Definition 4.1.1 and Definition 8.1.1]{TensorCategories}.
\end{remark}

\begin{lemma}
	The category $\overline{\Rep(\alpha_3)}$ is equivalent to the braided tensor category of super vector spaces.
	\begin{proof}
		Observe that $\End(\overline{\Phi^i})$ has a one-dimensional hom space. For $i = 1$ this is obvious. For $i = 2$, we note that the endomorphism $f$ is negligible. We also note that all homomorphisms between distinct $\Phi^i$ are negligible. 
		
		So, it is sufficient to prove that both categories have exactly $2$ non-isomorphic simple objects $O_0, O_1$ with one-dimensional Hom-spaces, such that $O_i \otimes O_j \cong O_{i+j}$, with indices modulo $2$, and that the braiding of those two categories coincide.
		For the super vector spaces, we note that this category is generated by a $1$-dimensional vector space $L_0$ and a $1$-dimensional vector space $L_1$ with $L_i \otimes L_j \cong L_{i +j},$
		and braiding given by
		\[ L_i \otimes L_j \ni a \otimes b \longmapsto (-1)^{ij} b \otimes a \in L_j \otimes L_i.\]
		
		For $\overline{\textbf{Rep}(\alpha_3)}$, we note that $\overline{\Phi^1}$ can play the role of $O_0$ and $\overline{\Phi^2}$ can play the role of $O_1$. The only nonobvious isomorphism is $O_1 \otimes O_1 \cong O_0$, which holds since
		\[ \Phi^2 \otimes \Phi^2 \cong \Phi^3 \oplus \Phi,\]
		in $\overline{\textbf{Rep}(\alpha_3)}$ using the decomposition of $\Phi^2 \otimes \Phi^2$ into symmetric and antisymmetric tensors.
		
		We remark that $\textbf{Rep}(\alpha_3)$ is equipped with the braiding $a \otimes b \mapsto b \otimes a$. This implies the braiding on $\overline{\Phi^2} \otimes \overline{\Phi^2}$ acts as multiplication by $-1$, as this corresponds to flipping the order in an antisymmetric tensor. On other tensor products of simple objects, the braiding acts trivially. This proves that we have the same braiding as for the category of super vector spaces.
	\end{proof}
\end{lemma}

\subsection{$J$-ternary algebras}

Wolfgang Hein \cite{HEIN} defined a $J$-ternary algebra as a module $M$ together with a trilinear mapping $(x,y,z) \mapsto xyz$ such that
\begin{align}
	xy(uvw) - uv(xyw) & = (xyu)vw + u(yxv)w, \\
	xyz - zyx &= zxy - xzy 
\end{align}
for all $x,y,z,u,v,w \in M$.
The $J$ stands for the linear span of the $M$-endomorphisms $\langle x, y \rangle$ that are defined as
\[ \langle x, y \rangle z = yzx - xzy.\]
This $J$ is a Jordan subalgebra of $\text{End}(M)$ with operation $x \cdot y = (xy + yx)/2$. Allison \cite{ALL76} defined $J$-ternary algebras in a similar fashion, but with a bigger emphasis on the Jordan part $J$. Any $J$-ternary algebra in the sense of Allison defines a $J$-ternary algebra in Hein's sense if one forgets the $J$. However, in Hein's setting, the Jordan algebra $J$ is uniquely determined by the triple product. Given that having different possible $J$ has no upside for what we want to do (as they will disappear entirely in the Lie superalgebra), we prefer to work with the definition of Hein.

We also assume that $J$ is a unital Jordan algebra. If we want to emphasize $J$, we write $(J,M)$ for the $J$-ternary algebra $M$ with Jordan algebra $J$.

\subsubsection{The associated Lie algebra}
Consider the left multiplication operator $L_{x,y} z = xyz$.
We define
\[ \text{InStr}(M) = \langle L_{x,y} \in \text{End}(M) | x,y \in M\rangle,\]
this is a Lie algebra acting on $M$. It also comes with an action on $J$ induced by $s \cdot \langle a , b \rangle = \langle s \cdot a,b \rangle + \langle a , s \cdot b \rangle$ for any $s \in \InStr(M)$ and all $\langle a , b \rangle \in J$.
Set \[\text{InDer}_J(M) = \{v \in \text{InStr}(M) | 0 = v \cdot 1 \},\]
using the unit $1$ of $J$.
The $\mathbb{Z}$-graded Lie algebra $L$ associated to $M$ is of the form
\[ L = J_- \oplus M_- \oplus \text{InStr}(M) \oplus M_+ \oplus J_+\]
with $J_\pm$ a $\pm 2$ graded copy of $J$, $M_\pm$ a $\pm 1$ graded copy of $M$, and with $\text{InStr}(M)$ a $0$-graded subalgebra of $L$. We will only need that $\InStr(M)$ acts naturally on $M_+$, that $[x,y] = L_{x,y}$ for $x \in M_+$ and $y \in M_-$, and that $[j_+,v_-] = (jv)_+$ for $j \in J$ and $v \in M$. For the precise algebra structure, see \cite{ALL76} or \cite{HEIN75}.

\begin{remark}
	The $J$-ternary algebras with \textit{unital} $J$ can be thought of as an algebraic structure out of which we can define a Lie algebra $L$ with a canonical inner derivation $f = \text{ad} \; 1_+$ such that $f^3 = 0$. Hence, if $\Phi$ is a field of characteristic $3$, we can interpret the associated Lie algebra as an element of $\textbf{Rep}(\alpha_3)$. In the next subsection, we will describe the Lie superalgebra one obtains after semisimplification in a bit more detail.
\end{remark}

\subsubsection{An interesting characterization}
An important characterization of the Lie algebras of $J$-ternary algebras was proved by Allison \cite{ALL76} over fields containing $1/6$. These Lie algebras contain a distinguished subalgebra $S = \langle 1_+, 1_-, [1_+,1_-] \rangle$ isomorphic to $\mathfrak{sl}_2(\Phi)$, where $1_\pm$ are the units of the Jordan algebra $J_{\pm 2}$, such that \[L = J \otimes (3) \oplus M \otimes (2) \oplus \text{InDer}_J(M) \otimes (1)\]
with $(i)$ the standard $i$-dimensional irreducible $\mathfrak{sl}_2$ module (and $J$, $M$, and $\InDer_J(M)$-trivial modules).
More precisely, Allison \cite[Theorem 1]{ALL76} proved that each Lie algebra $L$ with a Lie subalgebra $S \cong \mathfrak{sl}_2$ such that $L$ decomposes as
\[L = A \otimes (3) \oplus B \otimes (2) \oplus C \otimes (1)\]
for certain trivial $S$-modules $A, B$ and $C$ and the standard irreducible $i$-dimensional $S$-modules $(i)$, can be constructed from a $J$-ternary algebra if $C$ contains no ideals of $L$. The assumption that $C$ contains no ideals of $L$ can be removed, if we allow for central extensions of Lie algebras that can be constructed from $J$-ternary algebras \cite[Theorem 6.34]{Benkart2003LieAG}.

\subsubsection{Classification}
There exists an almost complete classification of the simple $J$-ternary algebras over algebraically closed fields depending on the degree of the associated Jordan algebra $J$.

\begin{definition}
	We call a $J$-ternary algebra $(J,M)$ \textit{simple} if the associated Lie algebra is simple. Hein proved this to be equivalent to the more natural notion of simplicity in \cite[Theorem 6.2]{HEIN81}. 
	Moreover, Hein \cite[Theorem 6.1]{HEIN81} proved that it is sufficient (and necessary\footnote{This definitely holds over fields of characteristic $0$. Hein's approach rests on being able to prove that $J$ is nondegenerate, and thus semisimple, in the terminology of \cite{Koe66}, to obtain that each ideal is derivation-invariant. This part of the proof is not entirely convincing.}) for $J$ to be simple and $\langle \cdot, \cdot \rangle$ to be nondegenerate.
\end{definition}

To be precise, the partial classification of $J$-ternary algebras (with simple $J$) that we will use can be stated as:
\begin{enumerate}
	\label{enumeration}
	\item Faulkner classified the Faulkner ternary algebras if $1/6 \in \Phi$, which are $J$-ternary algebras with $J$ equal to the basefield, i.e., a degree $1$ Jordan algebra; this classification was detailed and extended by Brown to characteristic $3$ \cite{Brown84} as indicated by \cite[Theorem 2.32]{ELD06};
	\item Hein classified $J$-ternary algebras of degree $\ge 3$ in \cite{HEIN81} if $1/2 \in \Phi$ as:
	\begin{enumerate}
		\item $J = H_n(C)$ with $C$ an associative composition algebra and $n \ge 3$, and $M$ given by $n \times m$ matrices over $C$;
	\end{enumerate}
	\item Hein classified $J$-ternary algebras of degree $2$ if $1/3 \in \Phi$ \cite{HEIN}, as:
	\begin{enumerate}
		\item $M$ is a tensor product of composition algebras;
		\item $M \cong A^n$ with $A$ a tensor product of associative composition algebras.
	\end{enumerate}
\end{enumerate}

\subsection{Semisimplification of the associated Lie algebra}

One can generalize the notion of a Lie algebra to an arbitrary braided tensor category. Normally, a Lie algebra has a bilinear, anticommutative Lie bracket $[\cdot,\cdot]$ that satisfies the Jacobi-identity $[a,[b,c]] + [b,[c,a]] + [c,[a,b]] = 0$. The bilinear Lie bracket, corresponds to a linear map $[\cdot,\cdot] : L \otimes L \longrightarrow L$. Using the braiding, we can generalize the anticommutativity and the Jacobi identity.

\begin{definition}
	Let $L$ be an object in a braided tensor category $\mathcal{C}$ with a homomorphism $[\cdot,\cdot]: L \otimes L \longrightarrow L$. This pair $(L,[\cdot,\cdot])$ is called an \textit{operadic Lie algebra} if the morphisms
	\[ [\cdot,\cdot] + [\cdot,\cdot] \circ b_{L,L}\]
	and 
	\[ [\cdot,[\cdot ,\cdot]](\text{Id} + (b_{L,L} \otimes 1)(1 \otimes b_{L,L}) + (1 \otimes b_{L,L})(b_{L,L} \otimes 1))\]
	are trivial, with $b$ the braiding on $\mathcal{C}$.
	We remark that operadic Lie algebras map to operadic Lie algebras under braided tensor functors. 
\end{definition}

\begin{lemma}
	If $(J,M)$ is a $J$-ternary algebra, then $\InDer_J(M) \oplus M$ with operation 
	\[ [(c,m), (d,n)] = ([c,d] + L_{m,n} + L_{n,m}, c(n)-  d(m))\]
	forms an operadic Lie superalgebra over fields $\Phi$ of characteristic $3$.
	\begin{proof}
		 Consider the $5$-graded Lie algebra \[L = J_- \oplus M_- \oplus \text{InStr}(M) \oplus M_+ \oplus J_+\] with $J_\pm$ and $M_\pm$ distinct copies of $J$ and $M$. 
		 We use the decomposition of the Lie algebra $L$ as a $\mathfrak{sl}_2$-module described earlier
		 \[L = J \otimes (3) \oplus M \otimes (2) \oplus (\text{Der}_J(M)) \otimes (1),\]
		 and observe that this becomes a $\Rep(\alpha_3)$-algebra with respect to $\ad 1_+$ for $1_+ \in J_+ \cap \mathfrak{sl}_2$.
		 
		 The semisimplification functor yields an operadic Lie superalgebra defined on $\text{InDer}_J(M) \oplus M$. The Lie bracket on $\text{InDer}_J(M)$ and the action of $\text{InDer}_J(M)$ on $M$ are preserved. For the bracket $M \otimes M \longrightarrow \text{Der}_J(M)$, we compute that
		 \[ [m \otimes \overline{\Phi^2}, n \otimes \overline{\Phi^2}] = [m_+,n_-] - [m_-,n_+] = [m_+,n_-] + [n_+,m_-], \]
		 with $m \otimes \overline{\Phi^2}$ denoting the embedding $\mu e_1 + \lambda e_2 \mapsto (\mu m_- + \lambda m_+ )$ of $\Phi^2$ in $L$. The Lie bracket $[m_+,n_-]$ is for a $J$-ternary algebra given by $L_{m,n}$.
	\end{proof}
\end{lemma}

Recall that $(J,M)$ is simple if and only if the Lie algebra is simple.

\begin{lemma}
	\label{lemma simple}
	Let $(J,M)$ be a $J$-ternary algebra.
	If $\InDer_J(M) \oplus M$ is simple, then 
	\[L = J_- \oplus M_- \oplus \InStr(M) \oplus M_+ \oplus J_+ \] is simple as well. Hence, $\InDer_J(M) \oplus M$ can only be simple if $(J,M)$ is simple as a $J$-ternary algebra.
	\begin{proof}
		We show that any nontrivial ideal of $L$ induces a nontrivial ideal of $\InDer_J(M) \oplus M$.
		
		Note that any ideal $I$ of $L$ has to be $\mathfrak{sl}_2$-invariant.
		Thus, any ideal $I$ decomposes as
		\[ I = (I \cap J_-) \oplus (I \cap M_-) \oplus (I \cap \InStr(M)) \oplus (I \cap M_+) \oplus (I \cap J_+),\]
		using the weight spaces of $d = \ad [1_+,1_-]$ and that the $-1$-weight space decomposes as 
		\[ J_+ \oplus M_- = \{l \in L : d(l) = - l, [1_+,l] = 0\} \oplus \{l \in L : d(l) = -l, [1_-,[1_+,l]] = l\},\]
		and using the analogous decomposition for the $+1$-weight space.
		Now, note that $J_- \oplus \InStr(M) \oplus J_+$ acts faithfully on $M_- \oplus M_+$.
		Hence, $I \cap (M_- \oplus M_+) = 0$ implies that $I$ is trivial.
		On the other hand $I \cap M_+ = M_+$ implies that $I = L$, as $L$ is generated by $M_+$ and $M_-$ and since $M_- = [1_-,M_+]$.
		
		If $I$ is nontrivial, set $I_+ = I \cap M_+$ and $D = I \cap \InDer_J(M)$. One verifies that $D \oplus I_+$ is a non-trivial ideal of $\InDer_J(M) \oplus M$.
%
%
	\end{proof}
\end{lemma}

\begin{remark}
	From here onwards, we will not write about general $J$-ternary algebras anymore.
	We will try to work with the Lie algebra instead.
	First, we sketch the result of applying semisimplification for the non-exceptional classes of $J$-ternary algebras.
\end{remark}

\subsection{Some classical Lie superalgebras from $J$-ternary algebras}

In this section, we will describe the Lie algebras related to $J$-ternary algebras of degree $\ge 3$ and the non-exceptional $J$-ternary algebras of degree $2$. We will not prove everything in detail, as none of the Lie algebras and Lie superalgebras are new. However, it is easy to verify that the Lie algebras we describe here correspond to the $J$-ternary algebras as classified in \cite{HEIN, HEIN81}.

The naming conventions used for the Lie superalgebras are the ones of Musson \cite{Musson}. For the superalgebra $\mathfrak{sl}(a, pk + a)$, which is not simple \cite[section 4]{Bou09}, we will still write $A(a,pk + a)$ for $\mathfrak{psl}(a,pk + a)$.

\subsubsection{Degree $\ge 3$}

We will look at $J$-ternary algebras of degree bigger than $3$ and show that these yield well-known superalgebras. These $J$-ternary algebras are given by $M_{r \times s}(\mathcal{C})$ with $r \ge 3$ and an associative composition algebra $\mathcal{C}$. We perform the construction starting from corresponding Lie algebra.
More precisely, consider \[M_{2r + s}(\mathcal{C}) = \begin{pmatrix}
	M_r(\mathcal{C}) & M_{r,s}(\mathcal{C}) & M_r(\mathcal{C}) \\
	M_{s,r}(\mathcal{C}) & M_{s}(\mathcal{C}) & M_{s,r}(\mathcal{C}) \\
	M_r(\mathcal{C}) & M_{r,s}(\mathcal{C}) & M_{r}(\mathcal{C}) \\
\end{pmatrix}\] for an associative composition algebra $\mathcal{C}$; we require $s$ to be even if $\mathcal{C}$ is $1$ dimensional. 
Take $J \in M_{s}(\mathcal{C})$ invertible such that $J^* = - J$, with $*$ denoting $A \mapsto \bar{A}^T$, and $G$ an invertible diagonal $r \times r$-matrix with entries in $\Phi$ (not $\mathcal{C}$).
Let 
\[ L = \left\{ a \in M_{2r + s}(\mathcal{C}) : a \begin{pmatrix}
	0 & 0 & - G^{-1} \\
	0 & J & 0 \\
	G^{-1} & 0 & 0 \\
\end{pmatrix} + \begin{pmatrix}
	0 & 0 & - G^{-1}\\
	0 & J & 0 \\
	G^{-1} & 0 & 0 \\
\end{pmatrix} a^* = 0 \right\}. \]
This is a Lie algebra of type $C_{(r + s/2)}, A_{2r + s}, D_{2r + s}$ depending on whether the dimension of $\mathcal{C}$ is $1,$ $2,$ or $4$ (for the type $D$, note that $L$ is preserving an anti-hermitian form over an algebra with symplectic involution). If $(2r + s + 1)$ is divisible by $3$ and $\mathcal{C} = E$, we should look at $[L,L]$ instead to obtain the Lie algebra generated by the $J$-ternary algebra.
The copies $M_1$ and $M_{-1}$ of $M$ in $L$ are formed by the elements 
\[ \begin{pmatrix}
	0 & x & 0 \\ 0 & 0 & J x^* G \\ 0 & 0 & 0
\end{pmatrix} \text{and} \begin{pmatrix}
0 & 0 & 0 \\ - Jy^* G & 0 & 0 \\ 0 & y & 0
\end{pmatrix}.\]
It is easy to verify that $[[x_1,y_{-1}],z_1] = xyz$ coincides with the product of \cite[Section 10]{HEIN} for $x,z \in M_1$ and $y \in M_{-1}$. 
Now, we semisimplify with respect to \[\begin{pmatrix}
	0 & 0 & I_r \\
	0 & 0 & 0 \\
	0 & 0 & 0
\end{pmatrix}.\] 

\begin{proposition}
	\label{prop: ABD}
	The Lie superalgebra one obtains is of type $A(r - 1,s - 1)$ if $\mathcal{C}$ is $2$-dimensional.
	If $\mathcal{C}$ is one dimensional one obtains $\mathfrak{osp}(r,s)$, which is of type $B$ or $D$ depending on whether $r$ is odd or even. If $\mathcal{C}$ is four dimensional, one obtains $\mathfrak{osp}(4s,4r)$, i.e., type $D(2s,2r)$. 
	\begin{proof}
		When $\mathcal{C}$ is two dimensional and $\Phi$ is algebraically closed, we know that $\mathcal{C} \cong \Phi \oplus \Phi$. In particular, if we project onto one of the copies of $\Phi$, we can see that we have a superalgebra generated by arbitrary $(r + s)^2$-block matrices of the form
		\[ \left( \begin{array}{ c | c}
			0   & B \\ \hline A & 0
		\end{array} \right),\]
		which means we get a superalgebra of type $A(r-1,s-1)$.
		
		In the other cases, we obtain a matrix that preserves the bilinear form $b$ induced by
		\[ \begin{pmatrix}
			- G^{-1} & \\ & J
		\end{pmatrix}\]
		since 
		
		\[ \left( \begin{array}{ c  c }
			0 & x \\  - Jx^*G & 0
		\end{array} \right)\begin{pmatrix}
			- G^{-1} & 0\\ 0& J
		\end{pmatrix} = \begin{pmatrix}
			0& x J \\ J x^* &0
		\end{pmatrix} = \begin{pmatrix}
			G^{-1} & 0 \\ 0& J
		\end{pmatrix} \begin{pmatrix}
			0 & x\\ - J x^*G & 0
		\end{pmatrix}^*.\]
		This is an orthosymplectic Lie algebra.
		Depending on the dimension of $\mathcal{C}$, the part corresponding to $G$ is the orthogonal part (1 dimensional) or the symplectic part (4 dimensional), cfr., \cite[Proposition (2.6)]{BookInvol}.	\end{proof}
\end{proposition}

\subsubsection{Non-exceptional degree $2$}

\label{subsubsec: nonexceptional degree 2}

For our convenience, we will describe the algebras one obtains over the algebraic closure. The types of superalgebras will remain the same if one works in full generality. As for the degree at least $3$ case, we immediately start from the Lie algebra.
We restrict ourselves to associative composition algebras $\mathcal{C}$. In sections \ref{sec: tpca} and \ref{sec: the con} we will look at what happens for non-associative composition algebras (which we consider the exceptional case). 

 So, let $A = \mathcal{C}_1 \otimes \mathcal{C}_2$ be a tensor product of associative composition algebras.
By considering traceless matrices $M$ in $M_{2+n}(A)$ such that
\[ M \begin{pmatrix}
	0 & 0 & 1 \\
	0 & I_n & 0 \\
	1 & 0 & 0
\end{pmatrix} + \begin{pmatrix}
	0 & 0 & 1 \\
	0 & I_n & 0 \\
	1 & 0 & 0
\end{pmatrix} M^*,\]
we obtain Lie algebras $L$ of the form $A_{n+1}, 2A_{n+1},C_{n+2}, A_{2n + 4}, D_{2n + 4}$ for $A = E, E \otimes E', Q, Q \otimes E, Q \otimes Q'$ with $E$ a composition algebras of degree $2$ and $Q$ of degree $4$ (writing $E'$ and $Q'$ for possible non-isomorphic composition algebras of the same degree); for the types $C$ and $D$ note that we preserve a hermitian form while $Q$ has a symplectic involution and $Q \otimes Q'$ has an orthogonal involution.
We remark that $E \otimes E'$ corresponds to a non simple Lie algebra. So, we will not consider $E \otimes E'$.

To obtain the relation with the $J$-ternary algebra, one can work as before with $[[x_1, (1 \otimes v \cdot y)_{-1}], z_1] = xyz$ to recover the triple product of \cite[Example preceding section $2$]{HEIN}. 

Now, we semisimplify with respect to an element
\[ \begin{pmatrix}
	0 & 0 & v \otimes 1 \\
	0 & 0 & 0 \\
	0 & 0 & 0 \\
\end{pmatrix}\]
with $v \otimes 1$ skew, i.e., $\overline{v \otimes 1} = - v\otimes 1.$ 

\begin{proposition}
	\label{prop: ACAD}
	The Lie superalgebra one obtains is of the form $A(0,n-1), C(n + 1), A(1,2n-1),$ or $D(2n,3)$ depending on whether $A = E,Q, E \otimes Q, Q \otimes Q'$.
	\begin{proof}
		To verify that the $A(0,n-1)$ and $A(1,2n-1)$ are correct, note that $E \cong \Phi \oplus \Phi$ whenever we work over an algebraically closed field. Hence, $M_{2 + n}(A) \cong M_{2 + n}(\Phi \oplus \Phi)$ or $M_{4 + 2n}(\Phi \oplus \Phi)$ depending on whether $A = E$ or $A = E \otimes Q$. Projecting onto one of the copies of $\Phi$, and applying the semisimplification, allows one to realize that we have a Lie superalgebra of the form $A(0,n-1)$ or $A(1,2n-1)$ depending on $A$, using a realization as $(1 + n)^2$-block-matrices or $(2 + 2n)^2$-block-matrices of the superalgebra.
		
		In the other cases, we obtain superalgebras formed by the $U$ in $M_{1 + n}(A)$ such that
		\[ U \begin{pmatrix}
			-v \\ & I_n
		\end{pmatrix} + \begin{pmatrix}
			-v \\ & I_n
		\end{pmatrix} U^* = 0,\]
	with the grading being the usual block-matrix-grading.
		Using that $A$ either has a symplectic or orthogonal involution, cfr., \cite[Definition (2.5)]{BookInvol}, we see that for $A = Q$ (symplectic, so that $M_n(A)$ with the hermitian transpose is symplectic) the block that only interacts with the $I_n$-part corresponds to a Lie algebra $\mathfrak{sp}_{2n}$, and that for $A = Q \otimes Q'$ (orthogonal) the $I_n$ part corresponds to $\mathfrak{o}_{4n}$. The Lie subalgebra corresponding to $v$ is either $\mathfrak{o}_2$ or $\mathfrak{sp}_6 \cong \mathfrak{o}_5$, which shall be clarified later in Lemma \ref{lemma orthogonal}. Upon closer examination, this matches precisely the classical block-matrix realization described in Musson \cite[section 2.3]{Musson}, assuming that one represents everything as matrices.	
	\end{proof}
\end{proposition}

\subsection{Degree $1$ and relation to \cite{ELD06}}

The semisimplification in degree $1$ was, in some sense, carried out by Elduque \cite{ELD06}.
In \cite[Theorem 2.32]{ELD06}, the classification of symplectic triple systems algebraically closed fields of characteristic $3$ was used. These symplectic triples systems are equivalent to Faulkner ternary algebras with nonzero alternating form \cite[Theorem 2.16]{ELD06}, which are precisely the degree $1$ $J$-ternary algebras.
For the symplectic triple systems, Elduque constructed a Lie superalgebra over fields of characteristic $3$ \cite[Theorem 3.1]{ELD06} that coincides precisely with the semisimplification we consider here.

One obtains many classical Lie superalgebras, described in the section before \cite[Theorem 3.2]{ELD06}. The new Lie superalgebras are listed in \cite[Theorem 3.2]{ELD06}. 

\begin{remark}
If one looks beyond algebraically closed fields, one can construct forms of the algebras described \cite[Theorem 3.2 (ii)-(v)]{ELD06} not considered in that article by making use of the structurable algebra associated with hermitian cubic norm structures, as introduced in \cite{Michiel}. In general, if one has a hermitian cubic norm structure $J$ over $K/R$, one can define an algebra $A = K \oplus J$ and look at $L_0 = \langle V_{x,sx} \rangle \subset \InStr(A)$ and define a Lie superalgebra structure on $L_0 \oplus A$.
\end{remark}

\begin{remark}
	In \cite{ELD06} one also looks at orthogonal triple systems. The Lie algebras obtained there in characteristic $3$ can also be thought of as coming from semisimplification applied to $\textbf{Rep}(\alpha_3) \boxtimes \mathbf{SVec} \longrightarrow \mathbf{SVec} \boxtimes \mathbf{SVec}$, using the Deligne tensor product $\boxtimes$, and then applying the obvious braided tensor functor to $\mathbf{SVec}$, as is obvious from \cite[Theorem 4.5]{ELD06} and \cite[Theorem 5.1]{ELD06}.
\end{remark}

\section{The tensor product of composition algebras}
\label{sec: tpca}

We will work with the tensor product of composition algebras by viewing them as structurable algebras. We proceed this way because the construction we want to carry out is somewhat easier to understand in terms of structurable algebras. At the end of this section, we will restate the general semisimplification idea in terms of this class of structurable algebras.

\subsection{The associated Clifford algebra and orthogonal Lie algebra}

Consider composition algebras $O_1, O_2$ over a field $\Phi$. Consider $A = O_1 \otimes O_2$ with an involution $a \otimes b \mapsto \bar{a} \otimes \bar{b}$. The subspace of skew elements $S$, i.e., $x \in O_1 \otimes O_2$ such that $\bar{x} = -x$ is given by $1 \otimes S_2 \oplus S_1 \otimes 1$ with $S_i$ the subspace of skew elements of $O_i$.

Consider the left multiplication operator $L_oa = oa$.
We use the norms $N_i$ of $O_i$, $i = 1,2$, to define a norm $N(a \otimes 1 + 1 \otimes b) = - N_1(a) + N_2(b)$ on $O_1 \otimes \Phi + \Phi \otimes O_2$. Let $N(s,v)$ be the polarization, i.e., $N(s,v) = N(s+v) - N(s) - N(v).$

The first thing we do is describe how we could interpret the skew elements of $A$ as a degree $2$ Jordan algebra.

\begin{lemma}
	\label{lem: Clifford rep}
	For all skew elements $s$ of $A$ and $w = 1 \otimes v$, we have
	$$ (L_sL_w)^2 + N(s,w)L_sL_w + N(s)N(w) = 0.$$
	\begin{proof}
		Let $s = a \otimes 1 + 1 \otimes b$.
		We compute
		\begin{align*}
			(L_sL_w)^2(y_1 \otimes y_2) & = a(ay_1) \otimes v(vy_2) + ay_1 \otimes b(v(vy_2)) \\  & \;\;\;+ ay_1 \otimes v(b(vy_2)) + y_1 \otimes b(v(b(vy_2))) \\
			& = N(a)N(v) y_1 \otimes y_2 - ay_1 \otimes N(v,b) vy_2 + y_1 \otimes b(v(b(vy_2)))\\ & \;\;\;+  y_1 \otimes b(b(v(vy_2))) - y_1 \otimes b(b(v(vy_2)))\\
			& = N(a)N(v) y_1 \otimes y_2 - ay_1 \otimes N(v,b) vy_2 - y_1 \otimes N(b,v) b(vy_2)\\ & \quad - N(b)N(v)y_1 \otimes y_2,\\
			& = - N(s)N(w)	y_1 \otimes y_2 - N(s,w) L_sL_v y_1 \otimes y_2,	
		\end{align*}
		where we used alternativity, and the fact that $x(\bar{y}z) + y(\bar{x}z) = N(x,y)z$, as this is a linearisation of $a(\bar{a}z) = N(a)z.$ 
	\end{proof}
\end{lemma}

We can generalize the previous lemma to also allow $w = v \otimes 1$ using $- N$.
One can prove the previous lemma even more generally. However, we will only be interested in $w = 1 \otimes v$ or $v \otimes 1$, since $L_wL_w$ is a multiple of the identity map for precisely those $w$.
Fix a skew element $w$ with $N(w)$ invertible.
We see that $L_wL_w/N(w) = \text{Id} \in \langle L_sL_w | s \text{ skew } \rangle = L_SL_w.$ 

\begin{definition}
	Consider a vector space $V$ and the tensor algebra $T(V)$ over $V$. Suppose that $q$ is a quadratic form on $V$.
	The \textit{Clifford algebra} is \[C(V) = T(V)/ (v \otimes v - q(v) : v \in V).\]If $V$ has dimension $n$, $C(V)$ has dimension $2^n$. Over algebraically closed fields $\Phi$, $C(V)$ is either isomorphic to $M_{2^{k}}(\Phi)$ or $M_{2^{l}}(\Phi \oplus \Phi)$ depending on whether $n = 2k$ or $n = 2l + 1$ \cite{Lam05}.
	Hence, if $V$ has dimension $2l + 1$ and we are working over the algebraic closure, $C(V)$ has a two non-isomorphic $2^{\ell}$-dimensional modules. Both modules are isomorphic as $\mathfrak{o}(q)$-modules and are denoted as $\spin_l$, using $\mathfrak{o}(q) \cong [V,V] \subset C(V)$ to define the action, as both are isomorphic the irreducible $2^\ell$-dimensional module of the central simple algebra $C_0(V)$ \cite[Chapter V, Theorem 2.4]{Lam05} (which contains $[V,V]$ as a Lie subalgebra).
\end{definition}

\begin{lemma}
	\label{lemma orthogonal}
	Consider the Clifford algebra $C(V)$ associated to the nondegenerate quadratic form $ s \mapsto N(s)N(w)$ defined on $V = \{ s \in S | N(s,w) = 0\}.$
	The special Jordan algebra $L_SL_w$ is isomorphic to the Jordan subalgebra $\Phi \oplus V$ of $C$ when we endow both with the product $a \cdot b = (ab + ba)/2$. Furthermore, $[L_SL_w,L_SL_w] \cong \mathfrak{o}(V)$.
	\begin{proof}
		First, note that $\Phi \oplus V$ is closed under the operation, as $v \cdot w = N(v,w)$ for $v,w \in V$ and $\phi \cdot v = \phi v$ for $\phi \in \Phi$ and $v \in \Phi \oplus V$.
		Secondly, observe that $(L_sL_w)^2 = -N(s)N(w)$ for $s$ such that $N(s,w) = 0$. So, we have two unital algebras defined on $\Phi 1 \oplus V$ where $1$ acts as a unit and where the multiplication on $V$ is determined by a symmetric bilinear map $V \times V \longrightarrow \Phi$, corresponding to $s \mapsto -N(s)N(w)$.
		
		We know that $[V,V] \cong \mathfrak{o}(V)$ in $C(V)$. Since $\Phi$ is central, this means that $[\Phi \oplus V,\Phi \oplus V] \cong \mathfrak{o}(V)$.
	\end{proof}
\end{lemma}

\subsection{The structurable algebras and associated Lie algebras}

The tensor product of composition algebras $A = O_1 \otimes O_2$ forms a \textit{structurable} algebra \cite[section 8.4]{ALL78} (or \cite[example 6.6]{ALL93} for fields of arbitrary characteristic), i.e., if $V_{x,y} z = (x\bar{y})z + (z\bar{y})x - (z\bar{x})y,$ then
$$ [V_{x,y},V_{u,v}] = V_{V_{x,y}u,v} - V_{u,V_{y,x}v}.$$

This structurable algebra corresponds to a $J$-ternary algebra with operation
\[ W_{x,y} = V_{x,vy},\]
for an element $v = s \otimes 1$ with $\bar{s} = -s$ and $N(s)$ invertible. We choose to work with the operator $V$ instead of $W$. Moreover, we prefer to work with the Lie algebra arising from $V$ instead of the isomorphic Lie algebra constructed from $W$. 

If $\Phi$ is a field of characteristic different from $3$, we can construct an associated Lie algebra \cite[section 3]{ALL79} of which we will recall the construction in just a moment. We can also work over fields of characteristic $3$, which follows from \cite[Section 5.4 and Lemma 3.16]{OKP}. The most important part is that the construction of the Lie algebra remains the same. 
This Lie algebra is isomorphic to the Lie algebra associated to $A$ as a $J$-ternary algebra. 

\begin{remark}
	\label{rmk: types Lie algebras}
	It is interesting to note that if $O_2$ is a Cayley algebra, and if $O_1$ is of dimension $2^i$, then we obtain a Lie algebra of type $F_4, E_6, E_7, E_8$ depending on whether $i = 0,1,2,3.$ 
	That this holds over fields of characteristic $0$, follows from \cite{ALL79}. In our arguments, we will see that over fields of characteristic $3$ the dimensions are bounded by the dimensions of $F_4$, $E_7$ and $E_8$, which is the only thing we need to make our argument work. 
	
	For $E_6$ we really need the type to be correct. However, Theorem \ref{thm: main} for $E_7$ or $E_8$ implies that the types for $E_7$ and $E_8$ are correct, since it lets us conclude that these are simple Lie algebras of the right dimension (and containing the right root system) using Lemma \ref{lemma simple}. Describing the Lie algebra of type $E_6$ then as a subquotient of $E_7$, allows us to see that the type for $E_6$ is also correct.
	
	To summarize, we will only work with the dimensions of the Lie algebras, except when working with $E_6$ where we also assume that the type is correct. Once we prove Theorem \ref{thm: main}, this stops being a problem.
\end{remark}

Set $D_0$ equal to the span of all $V_{x,y}$ with $x,y \in A$. We remark that there exist well-defined maps $\epsilon,\delta : D_0 \longrightarrow D_0$ such that 
\begin{enumerate}
	\item $V_{x,y}^\epsilon = - V_{y,x},$
	\item $ V_{x,y}^\delta (a\bar{b} - b\bar{a}) = (V_{x,y}a)\bar{b} - b(\overline{V_{x,y}a}) - a\overline{V_{x,y}b} + (V_{x,y}b)\bar{a},$
\end{enumerate}
 as proved\footnote{These equations also follow over fields of characteristic $3$ for $A$, as can be shown using \cite{OKP}. Namely, the definitions of $\epsilon$ and $\delta$ given in \cite{ALL79} are necessarily the correct ones to use if one wants to define an action of $D_0$ as inner derivations on the Lie algebra we shall define in just a moment. Over fields of characteristic $3$ this Lie algebra coincides with the Lie algebra of \cite[Lemma 3.16]{OKP}, which implies that the equations hold in general.} in \cite[section 1 and 2]{ALL79}.
The Lie algebra is defined on
\[ S_- \oplus A_- \oplus D_0 \oplus A_+ \oplus S_+,\]
with $A_\pm,S_\pm$ distinct copies of $A$ and $S$, which we see as a $\mathbb{Z}$-graded algebra.
The bracket of $D_0$ in this Lie algebra is given by $$d \mapsto (d^{\epsilon\delta}, d^\epsilon, \text{ad} \; d, d, d^\delta) : (S_- \times A_- \times D_0 \times A_+ \times S_+) \longrightarrow (S_- \times A_- \times D_0 \times A_+ \times S_+).$$
The other brackets are given by
\[ [x,y] = V_{x,y}, [s,t] = L_sL_t, [s,y] = sy, [t,x] = tx, [x,x'] = (x\bar{x}' - x'\bar{x})\]
for $x, x' \in A_+, y \in A_-, s \in S_+, t \in S_-$ where the result lies in the space which is uniquely determined by the $\mathbb{Z}$-grading. Using the Jacobi identity, one can write $L_sL_t$ as a linear combination of $V_{x,y}$.

There is an interesting subalgebra of $D_0$ given fixed invertible skew $v$, namely $ L_SL_v \oplus [L_SL_v,L_SL_v]$.

Now, we prove some properties about this Lie algebra that restate some basic properties about $J$-ternary algebras in terms of how we constructed the Lie algebra.

\begin{lemma}
	Consider a skew element $v$, then 
	$ [V_{x,vx},v] = 0$ holds for all $x \in A$.
	Furthermore, if $v$ is invertible, then
	$$ SV(v) = \langle V_{x,vy} - V_{y,vx} | x,y \in A\rangle = \langle L_sL_v | \bar{s} = -s \rangle $$
	maps bijectively to $S_-$ under
	$ d \mapsto [d,v].$
	Hence, for invertible skew $v$ we can decompose 
	$$ D_0 = \langle V_{x,vx} | x \in A \rangle \oplus SV(v).$$
	\begin{proof}
		We identify $v$ with an element of $S_-$. If $v$ is invertible, we identify its inverse with $v^{-1} \in S_+$.
		
		For $x,y,z \in A$, we compute that $[[v,[x,vy]],z] = [[vx,vy],z]$.
		Hence, $[V_{x,vy},v](z) = L_{[vx,vy]} (z)$.
		Substituting $y = x$ shows us that $[V_{x,vx},v] = 0$.
		
		So, we assume that $v$ is invertible.
		We compute that
		\[ [v^{-1},[v,[x,vy]]] = [v^{-1},[vx,vy]] = [x,vy] - [y,vx] = V_{x,vy} - V_{y,vy}.\]
		This means that $[v^{-1}/2,[v,d]] = d$ for $d \in SV(v)$.
		Hence, $\langle V_{x,vx} | x \in A\rangle  \oplus SV(v) \subseteq D_0$.		
		Each $V_{x,z}$ can be rewritten as $V_{x,v(l)}$ for some $l$, using the invertibility of $v$. And each $V_{x,vl}$ is an element of the direct sum, using $V_{x,vl} = (V_{x,vl} + V_{l,vx})/2 + (V_{x,vl} - V_{l,vx})/2$.
		
		We still need to prove that $SV(v) = \langle L_sL_v \rangle.$
		We compute that \[[v^{-1},[v,V_{x,vy}]] = [v,[{v^{-1}},V_{x,vy}]] = [v,[x,y]] = - L_{[x,y]}L_v,\]
		which shows that $2 L_{[x,y]}L_v = (V_{x,vy} - V_{y,vx}).$
	\end{proof}
\end{lemma}

The previous lemma shows that $SV(v)$ plays the role of the Jordan algebra $J$ for $J$-ternary algebras, while $\langle V_{x,vx} \rangle$ plays the role of $\text{InDer}_{J}(M) \subset \text{InStr}(M).$

\begin{lemma}
	Take a skew element $v$ in $A$ such that $v$ is invertible, and $a,b$ in $A$.
	Each $y \in A$ satisfies $[V_{y,vy},V_{a,vb}] = V_{V_{y,vy}a,vb} + V_{a,vV_{y,vy}b}$ and each element $d \in SV(v)$ satisfies $[d,V_{a,vb}] = V_{da,vb} - V_{a,v(db)}$.
	Furthermore, all $a,b,c,d \in A$ satisfy
	$$ [V_{a,vb},V_{c,vd}] = V_{V_{a,vb}c,vd} + V_{c,vV_{b,va}d}.$$
	\begin{proof}
		We compute
		$$ [e,[c,[v,d]]] = [[e,c],[v,d]] + [c,[[e,v],d]] + [c,[v,[e,d]]],$$
		for $e \in D_0$, $c,d \in A^+$ and $v \in S_-$ skew invertible.
		
		This proves the first part for $e = V_{x,vx}$ since $[e,v] = 0$.
		
		If $e = L_sL_v \in SV(v)$, then
		$e^\epsilon(vd) = - v(s(vd)) = - v(ed),$
		since $(L_sL_v)^\epsilon = L_vL_s$.
		So $[[e,v],d] = - 2v(ed)$.
		This proves the equation for $e \in SV(v)$.
		
		By using the direct sum $\langle V_{x,vy} : x,y \in A \rangle = SV(v) \oplus \{ V_{x,vx} : x \in A\}$ to rewrite $V_{a,vb}$, one proves the final equation.
	\end{proof}
\end{lemma}

\begin{lemma}
	\label{lemma 1}
	The space $\langle [L_sL_v,L_tL_v] : s,t \text{ skew}\rangle \subset \InStr(M) $ is contained in $\langle V_{x,vx} \rangle$.
	\begin{proof}
		By using $$[d,V_{a,vb}] = V_{da,vb} - V_{a,v(db)}$$
		for $d \in SV(v)$ and $a,b \in A$,
		one obtains
		$$ [e,V_{c,vd} - V_{d,vc}] = V_{ec,vd} + V_{d,v(ec)} - V_{c,v(ed)} - V_{ed,vc}.$$
		This proves the lemma, since $\langle L_sL_v \rangle = \langle V_{c,vd} - V_{d,vc} \rangle$.
	\end{proof}
\end{lemma}

\begin{remark}
	We restate the semisimplification idea from the previous section in this context.
	We have a Lie algebra
	\[ S_- \oplus A \oplus \text{InStr}(A) \oplus A \oplus S_+,\]
	with a decomposition of $\text{InStr}(A)$ as $S \oplus \langle V_{x,vx} \mid x \in A \rangle$ for arbitrary invertible $1 \otimes v \in S$.
	 With respect to this $v$ (now seen as an element of $S^+$), this Lie algebra is contained in $\text{Rep}(\alpha_3)$ and decomposes as
	 \[ (3) \otimes S \oplus (2) \otimes A \oplus (1) \otimes \langle V_{x,vx} \mid x \in A \rangle.\]
	 After semisimplification, we obtain an operadic Lie superalgebra
	 \[ L = \langle V_{x,vx} | x \in A \rangle \oplus A.\]
	 The bracket on the odd part is given by $(x,y) \mapsto V_{x,vy} + V_{y,vx}$.
\end{remark}

\section{The construction}
\label{sec: the con}

Let $A = O_1 \otimes O_2$ be the tensor product of composition algebras. Take $v \in O_1 \otimes O_2$ skew and invertible of the form $1 \otimes w$ or $w \otimes 1$.
Recall the action of the Clifford algebra $\text{Cl}(p)$, generated by $L_SL_v$ with $S$ the set of skew elements on $A$, for the quadratic form \[p: \{ s \in A | s \text{ skew, } s \perp v\} \longrightarrow \Phi\] given by $s \mapsto N(s)N(v).$
Recall that this induces an action of $\mathfrak{o}(p) \cong \langle [L_sL_v,L_tL_v] : s,t \in S\rangle$ on $A$.

If $O_1$ is a quaternion algebra, then $L_1 = S \cap O_1 \otimes 1$ forms a Lie subalgebra of $O_1 \otimes 1$. This Lie algebra acts on $A$ using right multiplication. The same can be done when $O_2$ is a quaternion algebra.
Similarly, if $O_1$ or $O_2$ is a composition algebra of degree $2$, then we consider $L_1 = \langle 3R_v \rangle$, i.e., the map $l \longrightarrow 3lv$, for skew invertible $v$. If $O_i \cong \Phi$ or $O_i$ is Cayley, we set $L_i \cong 0$. Thus, we define $R_L \cong L_1 \oplus L_2$ and note that it has dimension $6,4,3,2,1,$ or $0$ if the characteristic is different from $3$, while only $6,3,0$ are possible if the characteristic equals $3$.

\begin{lemma}
	\label{lem: description L0}
	For $A$ and $v$ as just introduced with $A$ not the tensor product of two composition algebras of degree $2$, we have
	$$ \langle V_{x,vx} \mid  x \in A \rangle \cong \mathfrak{o}(p) \oplus R_L.$$
	\begin{proof}
		First and foremost, note that we already proved that $\mathfrak{o}(p)$ is contained in the left hand side in Lemma \ref{lemma 1} and Lemma \ref{lemma orthogonal}.
		
		Secondly, $R_L \neq 0$ if either $O_1$ or $O_2$ is a composition algebra of degree $2$ or $4$. 
		Suppose that $O_1$ is a composition algebra of degree $4$.
		Take $w,s \in S \cap O_1 \otimes 1$ with $N(w,s) = 0$ and $N(s) \neq 0$. One computes that
		\[ (V_{s,ws} - N(s)V_{1,w})(y) = yw N(s).\]
		This proves that right multiplication by $w$ lies in $\langle V_{x,y} | x,y \in A \rangle$. One computes that any right multiplication $R_s$ by a skew element $s$ satisfies $R_s^\delta v = 0$, since $d^\delta(v) = \delta(v) + v\overline{d(1)}$. So, $R_w$ acts trivially on the $\pm 2$-graded parts of the Lie algebra associated with $A$. We know for sure that it lies in $\langle V_{x,vx} | x \in A \rangle$ since this last space is the kernel of $d \mapsto [d,v]$.
		If $O_1$ is a composition algebra of degree $2$, it is commutative and $w = v \otimes 1$. One computes that $V_{w,1}(y) = 3yw$, proving that $3R_w$ is contained in $D_0$. It also acts trivially on the $\pm 2$-graded part of the Lie algebra.
		This proves that $\mathfrak{o}(p) \oplus R_L$ is contained in $\langle V_{x,vx} | x \in A \rangle$.
		
		In characteristic $0$, one can compute \cite[8.4]{ALL79} the dimensions of the Lie algebra associated with $A$ using $\text{dim} \; L = 3 \cdot \text{dim} A + 2 \cdot \text{dim} S + l$ with $l$ being the dimension of the inner derivation Lie algebra of $O_1 \otimes 1 + 1 \otimes O_2$, which can be computed using the fact that this is a direct sum of the inner derivations of $O_1 \otimes 1$ and $1 \otimes O_2$, where each of these spaces is of dimension $0,0,3,14$ depending on the degree of the composition algebra\footnote{In the case of \cite[8.4]{ALL79} one assumes that one of the algebras is a Cayley algebra since the other algebras are not interesting. However, this does not change the formula for $l$.}.
		
		Moreover, if one looks at the simple analogs of these algebras over arbitrary characteristics, one knows that the center must lie in the part that acts trivially on the $\pm 2$-graded part, which we will prove to be the right multiplications using dimension arguments.
		Consider the table of dimensions in characteristic $0$ for the Lie algebra associated with $A$, where we use $E$ to denote an arbitrary composition algebra of degree $2$, $Q$ of degree $4$ and $O$ of degree $8$:
		\begin{center}
		\begin{tabular}{ccccc}
			Dim $L$ & $\Phi$ & $E$ & $Q$ & $O$ \\
			$E$ & 8 & X & 35 & 78\\
			$Q$ & 21 & 35 & 66 &  133 \\
			$O$ & 52 & 78 & 133 & 248 \\
		\end{tabular}	,
\end{center}	
	where we put an $X$ at $E \otimes E$ because this does not correspond to a simple Lie algebra, while all the others do.
	Furthermore, we do not have a row with $\Phi$ since the table is symmetric and $\Phi \otimes \Phi$ does not have any skew elements.
	If we go to arbitrary fields, these dimensions provide upper bounds, since the analogous construction over $\mathbb{Z}$ for the split composition algebras yields a free $\mathbb{Z}$-module. We remark that the dimensions will be preserved except for $8,35,78,$ as these correspond to Lie algebras of type $A_{3k + 2}$ or $E_6$ which have $1$-dimensional center over fields of characteristic $3$ (see, e.g., \cite{Bou09} for a table). For the ones of dimension $8$ and $35$, the link with Lie algebras of type $A_{3k + 2}$ was already clear in subsubsection \ref{subsubsec: nonexceptional degree 2}. For $E_6$, we recall Remark \ref{rmk: types Lie algebras} to indicate why we are justified in assuming the correct type.
	
	Computing the dimension \[\dim \langle V_{x,vx} | x \in A \rangle = \text{dim} L - 2 \text{dim} A - 3 \text{dim} S = \text{dim} \;\mathfrak{o}(\text{dim} S - 1) + \text{dim} R_L,\] proves the claim that only the right multiplications act trivially on the $\pm 2$-graded part (the first equality follows since $L = \langle V_{x,vx} \rangle \oplus A \otimes (2) \oplus S \otimes (3)$ as a $\mathfrak{sl}_2$-module). For the first row (and second column), note, firstly, that the types of the Lie algebras are $A_2$, $A_5$, and $E_6$, which have a one dimensional center over fields of characteristic $3$, and, secondly, that these are precisely the constructions involving the composition algebra $E$.
	We illustrate the computation for $O \otimes \Phi$:
	\[ 52 - 2 \cdot (\dim \Phi)(\dim O) - 3 \cdot (\dim \Phi + \dim O - 2) = 15 = (2 \cdot 3^2 - 3) = \dim \mathfrak{o}_6,\]
	and $O \otimes E$:
	\[ 78 - 2\cdot16- 3\cdot8 = 22 = (2 \cdot 3^2 + 3) + 1 = \dim \mathfrak{o}_7 + 1.\]
	This proves the direct sum decomposition claimed in the lemma.
	\end{proof}
\end{lemma}

In the lemma below, we characterize the bracket on the odd part of the algebra $A$ in the exceptional cases, i.e., we look at the map $(x,y) \longrightarrow V_{x,vy} + V_{y,vx}$, which is a $\InDer_J(M)$-invariant map over all fields.
When we used these characterizations, we will justify the modules we used to stand in for $A$. 

\begin{lemma}
	\label{lem: bracket char}
	Suppose that we work over $\mathbb{C}$ and let $S^2(V)$ be the symmetric square of a module $V$.
	\begin{itemize}
		\item There exists, up to rescaling, a unique $\mathfrak{o}_7$-invariant map \[S^2(\spin_7 \oplus \spin_7)  \longrightarrow \mathfrak{o}_7.\]
		\item There exists, up to rescaling, a unique $\mathfrak{o}_{13}$-invariant map \[S^2(\spin_{13}) \longrightarrow \mathfrak{o}_{13}.\]
		\item There exists, up to rescaling, a unique $\mathfrak{o}_{9}$-invariant map \[S^2(\spin_{9} \oplus \spin_{9}) \longrightarrow \mathfrak{o}_{9}.\]
		\item There exists, up to rescaling, a unique $(\mathfrak{o}_{9} \oplus \mathfrak{sl}_2)$-invariant map\[S^2(\spin_{9} \otimes (2)) \longrightarrow \mathfrak{sl}_{2}\] where $(2)$ denotes the standard $2$-dimensional irreducible $\mathfrak{sl}_2$-module.
	\end{itemize}
	\begin{proof}
		The first three statements follow immediately from character computations. Namely, in the first and third case, the Lie algebra does not occur as a quotient of $S^2(\spin_x)$ while it occurs exactly once as a summand of $\spin_x \otimes \spin_x$. In the second case, the Lie algebra occurs only once as a quotient of $S^2(\spin_{13})$. For the final case, remark that $(\mathfrak{o}_{9} \oplus \mathfrak{sl}_2)$-invariance implies that the map comes from a map $S^2(\spin_9) \longrightarrow \Phi$ and the unique map $S^2((2)) \longrightarrow \mathfrak{sl}_2$.
	\end{proof}
\end{lemma}

\begin{lemma}
	\label{lem: algebra defined over char 0}
	The product on $\mathfrak{g}(6,6) \cong \mathfrak{o}_{13} \oplus \spin_{13}$ can also be defined over $\mathbb{Z}[1/2]$ and the odd bracket corresponds to the unique $\mathfrak{o}_{13}$-invariant map 
		$S^2(\spin_{13}) \longrightarrow \mathfrak{o}_{13}$ over $\mathbb{C}$. 
	\begin{proof}
		This algebra $\mathfrak{g}(6,6)$ can be obtained using semisimplification as in \cite[Corollary 4.6.2]{KAN21}. We use that $e_1 + e_2$, as in the corollary, and the sum of the elements corresponding to opposite roots $e_{-1} + e_{-2}$, forms an $S$-triple satisfying the conditions of \cite[Theorem 1]{ALL76}. Hence, over $\mathbb{C}$, there exists a $J$-ternary algebra $(J,M)$ such that $\InDer_J(M) \oplus M$ with the (non-associative) product $(c,m) \cdot (d,n) = ([c,d] + L_{m,n} + L_{n,m}, c(n) - d(m))$ contains $\InDer_J(M)$ as a subalgebra and contains $M$ as a module for $\InDer_J(M)$. In addition, the product $M \otimes M \longrightarrow \InDer_J(M)$ is commutative and $\InDer_J(M)$-invariant. 		
		Using Lemma \ref{lem: bracket char}, we see that the operation on $M \otimes M$ is uniquely determined by being commutative and $\InDer_J(M)$-invariant. 
		
		The $J$-ternary algebra can be defined over $\mathbb{Z}$, using that the grading on the Lie algebra associated to the $J$-ternary algebra corresponds to a partition of the roots of $\mathfrak{e}_8$. If we have $1/2$ in our base ring we can split the $0$-graded part of the Lie algebra into $\InDer_J(M)$ and a copy of $J$ using $[e_{-2} + e_{-1}, [e_1 + e_2,d_j + d_i]] = 2 d_j$ for $d_j$ in the copy of $J$ and $d_i \in \InDer_J(M)$.
		Thus, we see that the semisimplification of \cite[Corollary 4.6.2]{KAN21} will yield the algebra $\InDer_J(M) \oplus M$.
	\end{proof}
\end{lemma}

\begin{theorem}
	Consider a field $\Phi$ of characteristic $3$ and let $A = O_1 \otimes O_2$ be a tensor product of composition algebras with the degrees of $O_1$ and $O_2$ not simultaneously equal to $2$.
	Then, $$ (\mathfrak{o}(p) \oplus R_L) \oplus A$$
	forms a Lie superalgebra $L$ under the described action and bracket $[x,y] = V_{x,vy} + V_{y,vx}$ on the odd part.
	\begin{proof} The only thing we have to prove is that these algebras satisfy the Jacobi identity in its homogeneous form. This is sufficient, as these algebras are necessarily the operadic Lie algebras obtained under the semisimplification functor.
		
		Now, we note that over algebraically closed fields and $O_1 \cong O_2$ an octonion algebra, we obtain the algebra $\mathfrak{g}(6,6) = \mathfrak{o}(13) \oplus \spin_{13}$ \cite[Table introduction]{Eld08}. To be precise, the even part was already described. The odd part has to be a single copy of the spin module using dimension considerations. Lemma \ref{lem: algebra defined over char 0} and the fact that our definition yields an algebra over arbitrary characteristics (for which only the Jacobi identity on the odd part might fail), show that we obtain $\mathfrak{g}(6,6)$.
		
		The other superalgebras arise as subquotients of $\mathfrak{g}(6,6)$ and are thus Lie superalgebras as well. Namely, any tensor product of composition algebras over the algebraic closure defines a subalgebra of the bioctonion algebra and thus a subalgebra of $\mathfrak{g}(6,6)$. The construction of the superalgebra we performed, agrees with this construction, if we make the even part act faithfully on the odd part. Hence, the Jacobi-identity in its homogeneous form follows from the same identity for $\mathfrak{g}(6,6)$.
			\end{proof}
\end{theorem}

Not every Lie superalgebra we constructed this way is simple. For example for $A = E \otimes \Phi$ we obtain a superalgebra of superdimension $(0,2)$ which cannot be simple.

\begin{remark}
	In the next theorem, we will establish that the Lie superalgebra we constructed is simple in certain cases. A table, summarizing which Lie superalgebras will not be simple (indicated with an $X$) and which types one obtains otherwise, is:
	\begin{center}
		\begin{tabular}{ccccc}
			Type/Simple & $\Phi$ & $E$ & $Q$ & $O$ \\
			$E$ & X & X & $\mathfrak{psl}(2,2)$ & $\mathfrak{g}(3,3)$\\
			$Q$ & $\mathfrak{osp}(2,2)$ & $\mathfrak{psl}(2,2)$ & $\mathfrak{osp}(4,4)$ &  $\mathfrak{el}(5,3)$ \\
			$O$ & $\mathfrak{psl}(4,1)$ & $\mathfrak{g}(3,3)$ & $\mathfrak{el}(5,3)$ & $\mathfrak{g}(6,6)$ \\
		\end{tabular}.
		\end{center}
\end{remark}

\begin{example}
	We already provided all the tools to prove the aforementioned isomorphisms. However, given the unique structure of $\mathfrak{el}(5,3)$ it might be useful to look a bit more carefully, as we need the precise module structure used in Lemma \ref{lem: bracket char} to be correct. Hence, the algebra one starts from is $A = O \otimes Q$ for an octonion algebra $O$ and a quaternion algebra $Q$. Thus, the odd part is $32$-dimensional.
	We already proved that the $0$-graded part is $\mathfrak{o}(p) \oplus R_Q$ for a $9$-dimensional quadratic form $p$ (since $10$ is the dimension of the skew part of $A$).
	Thus, whenever $Q$ is split, we can take idempotents $t$ and $1 - t$ in $Q$. Both $At$ and $A(1 - t)$ are $\mathfrak{o}(p)$-modules. Since $At$ is $16$-dimensional, it has to be the spin module $\spin_9$. It is not hard to verify that $A \cong \spin_9 \otimes (2)$ if one considers $Q$ as $2 \times 2$-matrices. Whenever $Q$ is split, we conclude
	\[ \InDer_J(M) \oplus A \cong (\mathfrak{o}(p) \oplus \mathfrak{sl}_2) \oplus (\spin_9 \otimes (2)).\]
	
	For $A = O \otimes E$, note that $A$ has to be a Clifford algebra module. Hence, $A$, which is $16$ dimensional, must decompose into $2$ copies of $\spin_7$ (as the quadratic form is $7$-dimensional). 
	\end{example}

\begin{theorem}
	\label{thm: main}
	Let $A$ be a tensor product of composition algebras.
	The Lie superalgebra $L$ associated to $A$ is simple if and only if
	$A$ is not isomorphic to $E \otimes E',$ $E \otimes \Phi$, or $\Phi \otimes \Phi,$ with $E$ and $E'$ composition algebras of degree $2$.
\begin{proof}
	For $E \otimes \Phi$, we obtain a Lie superalgebra of type $A(0,0)$, $\mathfrak{psl}(1,1)$ which is not simple, as observed in Proposition \ref{prop: ACAD}.
	For the other excluded possibilities, we did not even construct a Lie superalgebra.
	
	We now consider the remaining cases.
	We let $E, Q, O$ represent composition algebras of degree $2,4$ and $8$ respectively. For simplicity, we assume we are working over an algebraically closed field.
	\begin{enumerate}
		\item $A = Q \otimes \Phi$
		The Lie superalgebra obtained is isomorphic to $\mathfrak{osp}(2,2)$ (this is part of the broader class described in subsubsection \ref{subsubsec: nonexceptional degree 2}), which is simple.
		\item $A = Q \otimes E$
		We obtain $\mathfrak{o}_3 \oplus \mathfrak{sl}_2 \cong \mathfrak{sl}_2 \oplus \mathfrak{sl}_2$ as $0$-graded part. The full Lie superalgebra is isomorphic to $\mathfrak{psl}(2,2)$ (see subsubsection \ref{subsubsec: nonexceptional degree 2}), which is simple.
		\item $A = Q \otimes Q$
		We obtain a Lie superalgebra with $0$-graded part $$\mathfrak{o}(5) \oplus \mathfrak{sl}(2) \oplus \mathfrak{sl}(2) \cong \mathfrak{o}(5) \oplus \mathfrak{o}(4) \cong \mathfrak{sp}(4) \oplus \mathfrak{o}(4).$$ One checks that the Lie superalgebra is isomorphic to $\mathfrak{osp}(4,4)$ by comparing to subsubsection \ref{subsubsec: nonexceptional degree 2}.
		\item $O \otimes \Phi$
		We obtain a Lie superalgebra with $0$-graded part $\mathfrak{o}(6) \cong \mathfrak{sl}_4$ and the $1$-graded part is the simple module for the associated Clifford algebra. This module decomposes as the direct sum of two non-isomorphic modules for $\mathfrak{o}_6$. Using character computations one sees that the product on the odd part of the algebra is the unique (up to rescaling) $\mathfrak{sl}_4$-invariant commutative product over $\mathbb{C}$.
		This Lie superalgebra should thus be isomorphic to the simple algebra $\mathfrak{psl}(4,1)$.
		\item $E \otimes O, Q \otimes O, O \otimes O$
		By consulting the table in \cite{Eld08}, one can verify that the first and the third coincide with the Elduque Cunha superalgebras $\mathfrak{g}(3,3)$ and $\mathfrak{g}(6,6)$. The second one coincides with $\mathfrak{el}(5,3)$ \cite[12.2]{Bou09}. Those three are simple Lie superalgebras.
		 
		We already proved that this is the case for $\mathfrak{g}(6,6)$, using that we could define the product over $\mathbb{C}$. The semisimplifications of $\mathfrak{e}_6, \mathfrak{e}_7$, and $\mathfrak{e}_8$ performed in \cite{KAN21} with respect to $e_1 + e_2$ allow us to use the same arguments regarding definability over $\mathbb{Z}[1/2]$ as we did for $\mathfrak{e}_8$ and $\mathfrak{g}(6,6)$ and then use Lemma \ref{lem: bracket char}. \qedhere
	\end{enumerate}
\end{proof}
\end{theorem}

\bibliographystyle{alpha}
\bibliography{bib}
\end{document}